\documentclass[a4paper, oneside,11pt]{article}

\usepackage[latin1]{inputenc}
\usepackage[T1]{fontenc}
\usepackage[english]{babel}
\usepackage{amssymb}
\usepackage{amsmath}
\usepackage{amsthm}
\usepackage{graphicx}

\usepackage[all]{xy}
\usepackage{verbatim}
\usepackage{amsmath}
\usepackage{amsthm}
\usepackage{graphicx}
\usepackage{amssymb}
\usepackage{epstopdf}
\usepackage{url}
\usepackage{amsfonts}
\usepackage{todonotes}

\newcommand{\br}{\mathbb{R}}

\newcommand{\bt}{\mathbb{T}}

\newcommand{\e}{\varepsilon}

\newcommand{\bb}{\mathbb}

\newtheorem{thm}{Theorem}

\newtheorem{defi}[thm]{Definition}
\newtheorem{remark}[thm]{Remark}

\def\be{\begin{equation}}
\def\ee{\end{equation}}
\def\bea{\begin{eqnarray}}
\def\eea{\end{eqnarray}}

\usepackage[font=small,labelfont=bf, margin=1cm]{caption}

\usepackage{enumerate}

\usepackage{bbm}

\numberwithin{thm}{section}
\numberwithin{equation}{section}

\usepackage{tikz}
\usetikzlibrary{matrix}

\usepackage{hyperref}
\usepackage{mathtools}
\mathtoolsset{showonlyrefs}

\newcommand*\di{\mathop{}\!\mathrm{d}}

\newcommand{\RR}{\mathbb{R}}

\newcommand{\TT}{\mathbb{T}}

\newcommand{\ds}{\displaystyle}

\title{Recent developments on quasineutral limits for Vlasov-type equations}

\author{
Megan Griffin-Pickering 
  \thanks{Department of Mathematical Sciences, Durham University, Lower Mountjoy, Stockton Road, Durham DH1 3LE, UK. Email: \textsf{megan.k.griffin-pickering@durham.ac.uk}}
  \and
Mikaela Iacobelli
  \thanks{ETH Z\"urich, R\"amistrasse 101, 8092 Z\"urich, Switzerland. Email: \textsf{mikaela.iacobelli@math.ethz.ch}}
}

\begin{document}

\maketitle


\abstract{
Kinetic equations of Vlasov type are in widespread use as models in plasma physics.
A well known example is the Vlasov-Poisson system for collisionless, unmagnetised plasma.
In these notes, we discuss recent progress on the quasineutral limit in which the Debye length of the plasma tends to zero, an approximation widely assumed in applications.
The models formally obtained from Vlasov-Poisson systems in this limit can be seen as kinetic formulations of the Euler equations.
However, rigorous results on this limit typically require a structural or strong regularity condition.
Here we present recent results for a variant of the Vlasov-Poisson system, modelling ions in a regime of massless electrons.
We discuss the quasineutral limit from this system to the kinetic isothermal Euler system, in a setting with rough initial data.
Then, we consider the connection between the quasineutral limit and the problem of deriving these models from particle systems. We begin by presenting a recent result on the derivation of the Vlasov-Poisson system with massless electrons from a system of extended charges. 
Finally, we discuss a combined limit in which the kinetic isothermal Euler system is derived.
}


\section{Introduction}

Plasma is a state of matter consisting of an ionised gas, formed by the dissociation of a neutral gas under the influence of, for example, high temperatures or a strong magnetic field. Various mathematical models are available to describe plasma, corresponding to different physical regimes (such as typical length and time scales). Here we will focus on systems of Vlasov-Poisson type, which are kinetic equations describing dilute, collisionless, weakly magnetised plasmas.

The charged particles in a plasma typically fall into two distinguished types: electrons and positively charged ions.
The respective masses of these two species differ significantly -- note that the proton-to-electron mass ratio is of order $10^3$ \cite{Bellan}.
The result is a separation between the relevant timescales of evolution for the two species. As a consequence, it is a reasonable approximation to model the two species to some extent separately, and moreover the two species require different models. 

The best known version of the Vlasov-Poisson system is a kinetic model for the electrons in a plasma, evolving in a background of ions that are assumed to be stationary. This approximation is justified by the aforementioned separation of timescales.
For simplicity we leave aside the issue of boundary conditions by discussing the system posed on the $d$-dimensional flat torus $\bt^d$, which reads as follows:
\be \label{eq:VP}
(VP) : = \begin{cases}
\partial_t f + v \cdot \nabla_x f + E \cdot \nabla_v f = 0, \\
E = - \nabla_x U, \;
 - \Delta U = \rho_f - 1, \\ \ds
f \vert_{t=0} = f_0 , \int_{\TT^d \times \RR^d} f_0 (x,v) \di x \di v = 1.
\end{cases}
\ee

In these notes, we instead focus on a related model for the ions in a plasma. On the ions' timescale, the electrons are comparatively fast moving. In particular, the electron-electron collision frequency $\nu_e$ is much higher than the ion-ion collision frequency $\nu_i$. For example, Bellan \cite[Section 1.9]{Bellan} gives a relation of the form $\nu_{e} \sim (m_e/m_i)^{-1/2} \nu_i$ for plasmas with similar ion and electron temperatures, where $m_e$ and $m_i$ denote the masses of, respectively, a single electron and a single ion.
Thus, when the mass ratio $m_e/m_i$ is small, the frequency of electron-electron collisions can be significant even when ion-ion collisions are negligible.

In the \textbf{massless electrons} limit, the mass ratio $m_e/m_i$ is assumed to tend to zero, motivated by the fact that it is small in applications.
As a consequence, the electron collision frequency tends to infinity.
In the formal limiting regime, the electrons are thermalised,
instantaneously assuming their equilibrium distribution, which is a Maxwell-Boltzmann law of the form
\be
\rho_e \sim e^{q_e \beta_e \Phi} ,
\ee
where $q_e$ is the charge of a single electron, $\beta_e$ is the inverse electron temperature, and $\Phi$ is the ambient potential.

Combining the Vlasov-Poisson system \eqref{eq:VP} with a Maxwell-Boltzmann law for the electron distribution leads to the \textbf{Vlasov-Poisson system with massless electrons}, or VPME system. After an appropriate rescaling of physical constants, this reads as follows:
\be \label{eq:VPME}
(VPME) : = \begin{cases}
\partial_t f+ v \cdot \nabla_x f + E \cdot \nabla_v f = 0, \\
E = - \nabla_x U, \; \Delta U = e^U - \rho_f , \\ \ds
f \vert_{t=0} = f_0, \;  \int_{\TT^d \times \RR^d} f_0 (x,v) \di x \di v = 1.
\end{cases}
\ee
This model is used in the plasma physics literature to model ion plasma. For a more detailed introduction to the model in a physics context, see Gurevich and Pitaevsky \cite{Gurevich-Pitaevsky75}. The VPME system has been used to study the formation of ion-acoustic shocks \cite{Mason71, SCM}, the development of phase-space vortices behind these shocks \cite{BPLT1991}, and the expansion of plasma into vacuum \cite{Medvedev2011}, among other applications.

From a mathematical perspective, the VPME system has been studied less than the electron Vlasov-Poisson system \eqref{eq:VP}.
The systems differ through the additional exponential nonlinearity in the elliptic equation for the electrostatic potential in the VPME system. The nonlinearity of this coupling leads to additional difficulties.
For example, while the well-posedness theory of the Vlasov-Poisson system is well established (see for example \cite{Lions-Perthame, Loeper, Pfaffelmoser, Ukai-Okabe}), for the VPME system this theory was developed more recently. The existence of weak solutions was shown in $\RR^3$ by Bouchut \cite{Bouchut}, while global well-posedness was proved recently by the authors in \cite{IGP-WP}.

The massless electrons limit itself is not yet resolved in full generality. Bouchut and Dolbeault \cite{Bouchut-Dolbeault95} considered the problem for a one species model described by the Vlasov-Poisson-Fokker-Planck system.
Bardos, Golse, Nguyen and Sentis \cite{BGNS18} studied a two-species model represented by a system of coupled kinetic equations.
Under the assumption that this system has sufficiently regular solutions, in the massless electron limit they derive the Maxwell-Boltzmann law for the electron distribution, and a limiting system for the ions that is very similar to the VPME system \eqref{eq:VPME}, but with a time-dependent electron temperature.
We also refer to Herda \cite{Herda16} for the massless electron limit in the case with an external magnetic field.

In these notes, we summarise some recent progress on two problems related to the VPME system. In Section~\ref{sec:QN}, we consider the quasineutral limit, in which a characteristic parameter of the plasma known as the Debye length tends to zero. The limit of the VPME system in this regime is a singular Vlasov equation known as the kinetic isothermal Euler system. 
In Section~\ref{sec:particles} we consider the derivation of the VPME and kinetic isothermal Euler systems from a particle system. The underlying microscopic system consists of `ions', here represented as extended charges, interacting with each other and a background of thermalised electrons.

\section{Quasineutrality} \label{sec:QN}

\subsection{The Debye Length}

Plasmas have several important characteristic scales, one of which is the \textbf{Debye (screening) length}, $\lambda_D$. The Debye length has a key role in describing the physics of plasmas: broadly speaking, it governs the scale of electrostatic phenomena in the plasma. For example, it characterises charge separation within the plasma, describing the scale at which it can be observed that the plasma contains areas with a net positive or negative charge, and so is not microscopically neutral.

In terms of the physical constants of the plasma, the electron Debye length $\lambda_D$ is defined by
\be \label{def:Debye}
\lambda_D : = \left ( \frac{\epsilon_0 k_B T_e}{n_e q_e^2} \right )^{1/2}.
\ee
In the above formula, $\epsilon_0$ denotes the vacuum permittivity, $k_B$ is the Boltzmann constant, $T_e$ is the electron temperature and $n_e$ is the electron density.
The ions similarly have an associated Debye length, which may differ from the electron Debye length. It is defined by the formula \eqref{def:Debye}, replacing the electron density, temperature and charge with the corresponding values for the ions.

Since the Debye length is related to observable quantities such as the density and temperature, it can be found for a real plasma. Typically, $\lambda_D$ is much smaller than the typical length scale of observation $L$. The parameter
$\e := \lambda_D/L$
is therefore expected to be small. In this case the plasma is called \textbf{quasineutral}: since the scale of charge separation is small, the plasma appears to be neutral at the scale of observation. Quasineutrality is a very common property of real plasmas - for example Chen \cite[Section 1.2]{Chen} includes quasineutrality as one of the key properties distinguishing plasmas from ionised gases more generally. 

The significance for Vlasov-Poisson systems becomes apparent after a rescaling.
When written in appropriate dimensionless variables, the Vlasov-Poisson systems acquire a scaling of $\e^2$ in front of the Laplacian in the Poisson equation for the electric field. For example, the VPME system \eqref{eq:VPME} takes the form
\be \label{eq:VPME-quasi}
(VPME)_\e : = \begin{cases}
\partial_t f_\e + v \cdot \nabla_x f_\e + E \cdot \nabla_v f_\e = 0, \\
E = - \nabla_x U, \\
\e^2 \Delta U = e^U - \rho_{f_\e} , \\ \ds
f_\e \vert_{t=0} = f_\e(0), \;  \int_{\TT^d \times \RR^d} f_\e (0,x,v) \di x \di v = 1.
\end{cases}
\ee

In plasma physics literature, the approximation that $\e \approx 0$ is widely used. For this reason, it is important to understand what happens to the Vlasov-Poisson system in the limit as $\e$ tends to zero. This is known as the \textbf{quasineutral limit}. Taking this limit leads to other models for plasma known as kinetic Euler systems.

\subsection{Kinetic Euler Systems} \label{sec:KIE}

Formally setting $\e=0$ in the system \eqref{eq:VPME-quasi} results in the \textit{kinetic isothermal Euler system} (KIsE):
\be \label{eq:KE-iso}
(KIsE) :=
\begin{cases}
\partial_t f + v \cdot \nabla_x f - \nabla_x U \cdot \nabla_v f = 0, \\
U = \log{\rho_f}, \\ \ds
f \vert_{t=0} = f_0, \, \, \,  \int_{\TT^d \times \RR^d} f_0 (x,v) \di x \di v = 1.
\end{cases}
\ee
This system was described and studied in a physics context in \cite{GPP, GPP2, Gurevich-Pitaevsky75}. The name arises from the fact that, for monokinetic solutions $f$, of the form
\be \label{monokinetic}
f(t,x,v) = \rho(t,x) \delta_0(v - u(t,x))
\ee
for some density $\rho$ and velocity field $u$, the KIsE system is equivalent to the following isothermal Euler system:
\be \label{eq:Euler-iso}
(IsE) :=
\begin{cases}
\partial_t \rho + \nabla_x \cdot \left ( \rho u \right ) = 0, \\
\partial_t \left ( \rho u \right ) + \nabla_x \cdot \left (\rho u \otimes u \right ) - \nabla_x \rho = 0 .
\end{cases}
\ee

The KIsE system \eqref{eq:KE-iso} can be thought of as a kinetic formulation of the isothermal Euler system \eqref{eq:Euler-iso}
To see this, consider a solution in the form of a superposition of monokinetic profiles: let
\be \label{def:multi-fluid}
f(t,x,v) = \int_{\Theta} \rho_\theta (t,x) \delta_0(v - u_\theta (t,x)) \pi(\di \theta),
\ee
for a measure space $(\Theta, \pi)$ and a family of fluids $(\rho_\theta, u_\theta)_{\theta \in \Theta}$. 
The multi-fluid representation \eqref{def:multi-fluid} can be used in the case where $f$ has a density with respect to Lebesgue measure on $\TT^d \times \RR^d$. However, it can also accommodate more singular situations. For example, if $\pi$ is a sum of $N$ Dirac masses, then the distribution \eqref{def:multi-fluid} can be used to describe a system of $N$ phases.

With this multi-fluid representation in mind, consider the following system of PDEs for the unknowns $(\rho_\theta, u_\theta)_{\theta \in \Theta}$:
\be \label{eq:KIsE-mf}
(KIsE)_{MF} :=
\begin{cases}
\partial_t \rho_\theta + \nabla_x \cdot \left ( \rho_\theta u_\theta \right ) = 0, \\
\partial_t \left ( \rho_\theta u_\theta \right ) + \nabla_x \cdot \left (\rho_\theta u_\theta \otimes u_\theta \right ) = - \rho_\theta \nabla_x U, \\ \ds
U = \log \int_{\Theta} \rho_\theta (t,x) \pi(\di \theta) .
\end{cases}
\ee
Given a (distributional) solution of this multi-fluid system, the formula \eqref{def:multi-fluid} then defines a distributional solution of the KIsE system \eqref{eq:KE-iso}. Thus \eqref{eq:KIsE-mf} is a multi-fluid formulation of KIsE \eqref{eq:KE-iso} and KIsE is a kinetic formulation of the isothermal Euler system \eqref{eq:Euler-iso}.
The use of multi-fluid representations of this type for Vlasov-type equations is discussed, for example, in \cite{Zakharov, Grenier96, Brenier1999}. 

A system closely related to the KIsE system can be formally obtained by linearising the coupling $U = \log{\rho_f}$ between $U$ and $\rho_f$ around the constant density $1$: since $\log t\approx t-1$ for $t$ close to one, one gets
\be \label{eq:VDB}
(VDB) : = \begin{cases}
\partial_t f + v \cdot \nabla_x f - \nabla_x U \cdot \nabla_v f = 0, \\
U = \rho_f - 1 \\ \ds
f \vert_{t=0} = f_0, \, \, \,  \int_{\TT^d \times \RR^d} f_0 (x,v) \di x \di v = 1.
\end{cases}
\ee
This system was named the \textbf{Vlasov-Dirac-Benney} (VDB) system by Bardos \cite{Bardos}.
The name `Benney' was chosen due to a connection with the Benney equations for water waves, in particular as formulated by Zakharov \cite{Zakharov}.

The VDB system formally
has the structure of a general Vlasov equation, in which the potential $U$ is of the form $U = \Phi \ast_x (\rho_f-1)$ for some kernel $\Phi$. In this case, the kernel would be a Dirac mass; this is the origin of the reference to Dirac.
In particular, this demonstrates the additional singularity of the VDB system in comparison to the Vlasov-Poisson system: in the Vlasov-Poisson system the potential $U$ gains two derivatives compared to the density $\rho_f$, while in the VDB system this regularisation does not occur.

For the Vlasov-Poisson system for electrons \eqref{eq:VP}, the quasineutral limit leads to the following \textit{kinetic incompressible Euler} system (KInE):
\be \label{eq:KE-inc}
(KInE) :=
\begin{cases}
\partial_t f + v \cdot \nabla_x f - \nabla_x U \cdot \nabla_v f = 0, \\
\rho_f = 1, \\ \ds
f \vert_{t=0} = f_0, \, \, \,  \int_{\TT^d \times \RR^d} f_0 (x,v) \di x \di v = 1.
\end{cases}
\ee
The force $- \nabla_x U$ is defined implicitly through the incompressibility constraint $\rho_f = 1$, and may be thought of as a Lagrange multiplier associated to this constraint.
The system \eqref{eq:KE-inc} was discussed by Brenier in \cite{Brenier1989} as a kinetic formulation of the incompressible Euler equations.

All three kinetic Euler systems described above \eqref{eq:KIsE-mf}, \eqref{eq:VDB}, \eqref{eq:KE-inc} as well as the two Vlasov-Poisson systems \eqref{eq:VP},\eqref{eq:VPME}, have a large family of stationary solutions: the spatially homogeneous profiles $f(t,x,v) = \mu(v)$. As is well-known for the Vlasov-Poisson system, some of these profiles may be unstable \cite{Penrose}.
For the kinetic Euler systems, the corresponding linearised problems have unbounded unstable spectrum: see \cite{Bardos-Besse, Bardos-Nouri, Han-Kwan-Nguyen}. As a consequence, they are in general ill-posed. For example, ill-posedness in Sobolev spaces was shown for the VDB system by Bardos and Nouri \cite{Bardos-Nouri}. Han-Kwan and Nguyen \cite{Han-Kwan-Nguyen} further extended this by showing that the solution map cannot be H\"{o}lder continuous with respect to the initial datum in Sobolev spaces, for both the VDB system \eqref{eq:VDB} and the KInE system \eqref{eq:KE-inc}.
See also Baradat \cite{Baradat} for the generalisation when the unstable profile $\mu$ is only a measure.

Due to these instability properties, well-posedness results for the kinetic Euler systems typically involve either a strong regularity restriction or a structural condition. For instance, in the monokinetic case one may appeal to the results known for the corresponding Euler system.

Without imposing any structural condition, the most general results available are in analytic regularity. Local existence of analytic solutions for the VDB system was proven by Jabin and Nouri \cite{Jabin-Nouri} in the one-dimensional case, and also follows from \cite[Section 9]{Mouhot-Villani}.
Bossy, Fontbona, Jabin and Jabir \cite{BFJJ} proved an analogous result for a class of kinetic equations involving an incompressibility constraint, generalising the KInE system \eqref{eq:KE-inc} to include, for example, noise terms.
Local existence of analytic solutions for the multi-fluid system corresponding to KInE \eqref{eq:KE-inc} was shown by Grenier \cite{Grenier96} as part of a study of the quasineutral limit; note that, due to the multi-fluid formulation, the required regularity is only imposed in the $x$ variable.

In Sobolev regularity, local well-posedness is known for the VDB system for initial data satisfying a Penrose-style stability criterion, following the results of Bardos and Besse \cite{Bardos-Besse} and Han-Kwan and Rousset \cite{Han-Kwan-Rousset}. We do not know of any global-in-time existence results for any of the kinetic Euler systems \eqref{eq:KE-iso}, \eqref{eq:VDB} or \eqref{eq:KE-inc}.

The VDB system also appears in the semiclassical limit of an infinite dimensional system of coupled nonlinear Schr{\"{o}}dinger equations: for more details, see for example 
\cite{Bardos-Besse, Bardos-Besse2015, Bardos-BesseSC}. See also \cite{Carles-Nouri, Ferriere} for discussion of semiclassical limits involving the KIsE model.

\subsection{Failure of the Quasineutral Limit}

The mathematical justification of the quasineutral limit is a non-trivial problem, since in general the limit can be false.
The failure of the limit can be linked to 
known phenomena in plasma physics.
We note for instance the example of Medvedev \cite{Medvedev2011} regarding the expansion of ion plasma into vacuum. For a one-dimensional hydrodynamic model it is found that the quasineutral approximation $U = \log \rho$ is not valid everywhere, and this is corroborated by numerical simulations for a kinetic model. 

Another important issue, well-known in plasma physics, is the `two stream' instability. 
From a physics perspective, this instability is typically introduced through a model problem in which two jets of electrons are fired towards each other (whence the name). Configurations of this kind are known to be unstable (see for example \cite[Section 5.1]{Bellan}, \cite[Section 6.6]{Chen}), with the resulting dynamics producing a vortex-like behaviour in phase space. See \cite{BNR} for simulations and experimental results on this phenomenon.
The streaming instability is seen in kinetic models by considering profiles with a `double bump' structure in the velocity variable.
These profiles are unstable for the linearised problem in the Penrose sense discussed above.

The relevance of instability for the quasineutral limit can be indicated by looking at a time rescaling of the Vlasov-Poisson system. 
If $f$ is a solution of the unscaled Vlasov-Poisson system \eqref{eq:VP}, then $ f_\e(t,x,v) = f\left ( \frac{t}{\e}, \frac{x}{\e}, v \right )$ is a solution of the system with quasineutral scaling. The limit as $\e$ tends to zero is thus a form of long time limit.
Grenier outlined this obstruction to the quasineutral limit in \cite{Grenier96, Grenier99}, for a one-dimensional two-stream configuration. 
Subsequently, Han-Kwan and Hauray \cite{Han-Kwan-Hauray} constructed counterexamples to the quasineutral limit in the Sobolev spaces $H^s$ for arbitrary large $s$, by considering initial data around unstable profiles.

\subsection{Results on the Quasineutral Limit}

Positive results on the quasineutral limit can be categorised along the lines of the well-posedness results known for the kinetic Euler systems; these problems are closely related. 
The mathematical study of the quasineutral limit can be traced back to the 90s, with the works of Brenier and Grenier \cite{Brenier-Grenier94} and Grenier \cite{Grenier95}, using an approach based on defect measures, and the result of Grenier \cite{Grenier99} for the one-dimensional case.

A particular case is the `cold electrons' or `cold ions' regime, in which the initial data for the Vlasov-Poisson system is assumed to converge to a monokinetic profile. The limiting kinetic Euler system is therefore reduced to its corresponding Euler system.
Brenier \cite{Brenier2000} and Masmoudi \cite{Masmoudi2001} considered the electron case, from the Vlasov-Poisson system to the incompressible Euler equations. Han-Kwan \cite{Han-Kwan2011} considered the ions case, from the VPME system to the isothermal Euler equations.
See also the work of Golse and Saint-Raymond \cite{Golse-SR2003}, obtaining a `2.5 dimensional' incompressible Euler system through a combined quasineutral and gyrokinetic limit (a limit of strong magnetic field).

In \cite{Grenier96}, Grenier proved the quasineutral limit from the electron Vlasov-Poisson system to KInE in analytic regularity. The result is framed in terms of the corresponding multi-fluid formulations. If the initial data for the multi-fluid Vlasov-Poisson system are uniformly analytic in $x$, then the quasineutral limit to the multi-fluid KInE system holds locally in time.
By the same techniques, similar results can be shown for the ion quasineutral limits, obtaining the VDB and KIsE systems, as observed in \cite{IHK1}, in the discussion after Proposition 4.1.

Under a Penrose-type stability criterion, Han-Kwan and Rousset \cite{Han-Kwan-Rousset} proved that the quasineutral limit holds in Sobolev regularity, for the passage from a variant of the VPME system, with linearised Poisson-Boltzmann coupling for the electric field, to the VDB system.

\subsection{Quasineutral Limit with Rough Data} \label{sec:result-QN}

An alternative direction for relaxing the regularity constraint for the quasineutral limit was investigated in a series of works, by Han-Kwan and the second author \cite{IHK2,IHK1} and by the authors \cite{GPI20}.
In this setting, one considers rough initial data (measures in the one-dimensional case, $L^\infty$ for $d=2,3$) that are small perturbations of the uniformly analytic case. The smallness of the perturbation is measured in a Wasserstein (Monge-Kantorovich) distance.

\begin{defi}[Wasserstein Distances] \label{def:Wass}
Let $p \in [1, \infty)$. Let $\mu$ and $\nu$ be probability measures on $\TT^d \times \RR^d$ for which the moment of order $p$ is finite.
Then the $p$\textsuperscript{th} order Wasserstein distance between $\mu$ and $\nu$, $W_p(\mu,\nu)$, is defined by
\be \label{def:MKW}
W_p(\mu, \nu) = \left ( \inf \int_{(z_1,z_2) \in (\TT^d \times \RR^d)^2} d(z_1, z_2)^p \, \di \pi(z_1, z_2) \right )^{1/p}, 
\ee
with the infimum taken over measures $\pi$ on $(\TT^d \times \RR^d)^2$ such that for all Borel sets $A \subset \TT^d \times \RR^d$,
\be
\pi(A \times \TT^d \times \RR^d) = \mu(A), \qquad \pi( \TT^d \times \RR^d \times A ) = \nu(A),
\ee
and $d$ denotes the standard metric on $\TT^d \times \RR^d$.
\end{defi}

The article \cite{IHK1} deals with the one-dimensional case for both electron and ion models, while  
in higher dimensions $d=2,3$, the limit for the electron models is considered in \cite{IHK2}. 
Then, for the VPME system, we proved a rough data quasineutral limit in \cite{GPI20}.

Below we give the statement of this result. We use the notation $\mathbf{\overline{\exp}_n} $ to denote the $n$-fold iteration of the exponential function, for example
\be
\mathbf{\overline{\exp}_3} (x) = \exp \exp \exp (x) .
\ee
We also use the analytic norms $\lVert \cdot \rVert_{B_\delta}$, defined for $\delta > 1$ by
\be
\lVert g \rVert_{B_\delta} : = \sum_{k \in \bb{Z}^d} |\hat g(k)| \delta^{|k|} ,
\ee 
where $\hat g(k)$ denotes the Fourier coefficient of $g$ of index $k$.

\begin{thm}[Quasineutral limit] \label{thm:quasi-summary}
Let $d = 2, 3$. Consider initial data $f_\e(0)$ satisfying the following conditions:
\begin{itemize}
\item (Uniform bounds) $f_\e(0)$ is bounded and has bounded energy, uniformly with respect to $\e$: for some constant $C_0>0$,
\be \label{unif-energy}
\lVert f_\e(0) \rVert_{L^{\infty}(\bt^d \times \br^d)}  \leq C_0, \qquad
\frac{1}{2}\int_{\TT^d \times \RR^d} |v|^2 f \di x \di v + \frac{\e^2}{2} \int_{\TT^d} |\nabla U |^2 \di x +  \int_{\TT^d} U e^{U} \di x \leq C_0 .
\ee

\item (Control of support) There exists $C_1>0$ such that 
\be \label{quasi:data-spt}
f_\e(0, x, v) = 0 \qquad \text{for } \; |v| > \exp(C_1 \e^{-2}) .
\ee
\item (Perturbation of an analytic function) There exist $g_\e(0)$ satisfying, for some $\delta > 1$, $\eta>0$, and $C>0$,
\be \label{analytic-assumptions}
\sup_{\e > 0} \sup_{v \in \br^d} (1 + |v|^{d+1}) \lVert g_\e(0, \cdot, v) \rVert_{B_\delta} \leq C , \qquad \sup_{\e > 0} \left \| \int_{\br^d} g_\e(0, \cdot, v) \di v - 1  \right \|_{B_\delta} \leq \eta ,
\ee
as well as the support condition \eqref{quasi:data-spt}, such that, for all $\e > 0$,
\be \label{Wass-rate}
W_2(f_\e(0), g_\e(0)) \leq \left [ \,\mathbf{\overline{\exp}_4}  (C \e^{-2}) \right ]^{-1} 
\ee
for $C$ sufficiently large with respect to $C_0, C_1$.
\item(Convergence of data) $g_\e(0)$ has a limit $g(0)$ in the sense of distributions as $\e \to 0$.
\end{itemize}
Let $f_\e$ denote the unique solution of \eqref{eq:VPME-quasi} with bounded density and initial datum $f_\e(0)$. Then there exists a time horizon $T_* > 0$, independent of $\e$ but depending on the collection $\{ g_{0,\e} \}_\e$, and a solution $g$ of \eqref{eq:KE-iso} on the time interval $[0, T_*]$ with initial datum $g(0)$, such that
\be
\lim_{\e \to 0}\, \sup_{t \in [0, T_*]} W_1(f_\e(t), g(t)) = 0 .
\ee
\end{thm}

\begin{remark}
As an example of a choice of initial data satisfying these assumptions, consider any compactly supported, spatially homogeneous profile $\mu = \mu(v) \geq 0$ with unit mass. Then
\be
f_\e(0) = \mu(v) \left (1 + \sin(2 \pi N_\e x_1) \right ), \quad N_\e\gtrsim \mathbf{\overline{\exp}_4}(C \e^{-2})
\ee
satisfies the assumptions of Theorem~\ref{thm:quasi-summary}.
\end{remark}

\subsection{Remarks on the Strategy} \label{sec:strategy-QN}

The strategy of proof for the rough data quasineutral limits \cite{IHK1, IHK2, GPI20} is based on stability results for the Vlasov-Poisson systems in Wasserstein distances.
Stability results of this type have been known for Vlasov-type equations since the work of Dobrushin \cite{Dobrushin} for the case of Lipschitz force kernels.

The Vlasov-Poisson case was considered by Loeper \cite{Loeper}, for solutions whose mass density $\rho_f$ is bounded in $L^\infty$.
This is an estimate of the form
\be
W_2(f_1(t), f_2(t)) \leq \mathcal{F} \left [ W_2(f_1(0), f_2(0)), \;\max_{i=1,2} \| \rho_{f_i} \|_{L^\infty([0,t] \times \TT^d)} \right ], 
\ee
for some suitable $\mathcal{F}.$
The corresponding estimate for the VPME system was proved recently in \cite{IGP-WP}.

The proof of Theorem~\ref{thm:quasi-summary} relies on a quantification of the $W_2$ stability estimate in terms of $\e$.
This has two steps: first, the stability estimate itself is quantified, in the sense that
\be
W_2(f_\e^{(1)}(t), f_\e^{(2)}(t)) \leq \mathcal{F}_\e \left [ W_2(f_\e^{(1)}(0), f_\e^{(2)}(0)), \;\max_{i=1,2} \| \rho_{f_\e^{(i)}} \|_{L^\infty([0,t] \times \TT^d)} \right ] .
\ee
Then, a bound is proved for the mass density $\| \rho_{f_\e^{(i)}} \|_{L^\infty([0,t] \times \TT^d)} $ in terms of the initial data. This is achieved by controlling the rate of growth of the support of a solution $f_\e$ in terms of the initial data, via an analysis of the characteristic trajectories of the system. This is the reason for the compact support assumption in Theorem~\ref{thm:quasi-summary}.

The quantified stability estimate is then used to make a perturbation around the analytic regime.
More specifically, we consider the analytic functions $g_\e(0)$ defined in the statement as initial data for the VPME system \eqref{eq:VPME}.
The assumptions \eqref{analytic-assumptions} are chosen precisely so that the resulting solutions $g_\e$ satisfy the quasineutral limit: on some time interval $[0, T_*]$, as $\e$ tends to zero, $g_\e$ converges to a solution $g$ of the KIsE system \eqref{eq:KE-iso}.
This follows from the techniques of Grenier \cite{Grenier96}, and implies convergence in a Wasserstein distance.

The proof is concluded by the triangle inequality:
\be
W_1(f_\e(t), g(t)) \leq W_1(f_\e(t), g_\e(t)) + W_1(g_\e(t), g(t)),
\ee
choosing the envelope of initial data \eqref{Wass-rate} so that the perturbation term $W_1(f_\e(t), g_\e(t)) $ vanishes in the limit.

\section{Derivations from Particle Systems} \label{sec:particles}

It is a fundamental problem to derive effective equations, such as Vlasov-Poisson systems, from the physical systems they are intended to describe.
In a reasonably general setting, we may consider a system of $N$ point particles with binary interactions.
The dynamics of such a system are modelled in classical mechanics by a system of ODEs of the following form, describing the phase space positions $(X_i, V_i)_{i=1}^N$ of the particles: 
\be \label{ODE-gen}
\begin{cases}
\dot X_i = V_i \\ \ds
\dot V_i = \alpha(N) \sum_{j \neq i} \nabla W(X_i - X_j) + \nabla V(X_i) .
\end{cases}
\ee
In this setting $\nabla W$ denotes the interaction force between pairs of particles, which here depends only on the spatial separation of the particles and is derived from an interaction potential $W$. We also include an external force $\nabla V$.
The parameter $\alpha(N)$ rescales the system with $N$ and can be thought of as a rescaling of the physical constants of the system. The choice of $\alpha(N)$ determines the model that is obtained as $N$ tends to infinity.

The case $\alpha(N) = 1/N$ is known as the \textbf{mean field limit}.
The formal limiting system is the Vlasov-type equation
\be \label{eq:vlasov-MFL}
\partial_t f + v \cdot \nabla_x f + (\nabla W \ast_x \rho_f + \nabla V) \cdot \nabla_v f = 0 ,
\ee
in the sense that the empirical measures $\mu^N$ defined by the formula
\be \label{def:mu}
\mu^N : = \frac{1}{N} \sum_{i=1}^N \delta_{(X_i, V_i)} 
\ee
are expected to converge to a solution of the Vlasov equation \eqref{eq:vlasov-MFL} in the limit as $N$ tends to infinity.
The Vlasov-Poisson system fits into this framework by choosing $\nabla V = 0$ and $\nabla W$ to be the Coulomb kernel $K$ on the torus $\TT^d$. This is the function $K = - \nabla G$, where $G$ satisfies
\be
- \Delta G = \delta_0 - 1 \qquad \text{on } \TT^d.
\ee
The corresponding microscopic system \eqref{ODE-gen} then describes a system of interacting electrons modelled as point charges, while \eqref{eq:vlasov-MFL} is the Vlasov-Poisson system \eqref{eq:VP}.

To derive the VPME system, a natural choice for the underlying microscopic system is to consider the dynamics of $N$ ions, modelled as point charges, in a background of thermalised electrons.
On the torus, this is modelled by an ODE system of the form
\be \label{eq:ODE-VPME}
 \left \{
\begin{array}{l}
\dot{X}_i = V_i \\ \ds
\dot{V}_i = \frac{1}{N} \sum_{j \neq i}^N K (X_i - X_j) - K \ast e^U  ,
\end{array}
\right.
\ee
where the electrostatic potential $U$ satisfies
\be \label{ODE-U}
\Delta U = e^{U} - \frac{1}{N} \sum_{i=1}^N \delta_{X_i} .
\ee
We can think of this system as being of the form \eqref{ODE-gen} by taking $\nabla W = K$ and an `external' force $\nabla V = K \ast e^U$, even though $\nabla V$ is not truly external due to its nonlinear dependence on the particle configuration through $U$. 
In this way it can be seen that the VPME system formally describes the limit as $N$ tends to infinity.

Other choices are possible for $\alpha(N)$, in which case the limit as $N$ tends to infinity may produce models of other forms. 
This approach can be used to derive the kinetic Euler systems discussed above in Subsection~\ref{sec:KIE}.
In the papers \cite{IGP1, GPI20}, the scaling $\alpha(N)\approx \frac{1}{N \log \log N}$ is used to derive the kinetic Euler systems \eqref{eq:KE-inc} and \eqref{eq:KE-iso}. The method is based on passing via the associated Vlasov-Poisson system, and this limit can thus be thought of as a simultaneous mean field and quasineutral limit.
In the recent paper  \cite{IHK-SCMFL}, a similar limit is proved in the monokinetic regime, to derive the incompressible Euler equations.

\subsection{Mean Field Limits}

For a detailed survey of mathematical results on the mean field limit, see \cite{Golse, Jabin_MFLreview}.
For our purposes we emphasise that the theory of mean field limits depends on the regularity of the interaction force $\nabla W$ chosen in the system \eqref{eq:vlasov-MFL}.

Early contributions on the problem include the works of Braun-Hepp \cite{Braun-Hepp}, Neunzert-Wick \cite{Neunzert-Wick} and Dobrushin \cite{Dobrushin}. In particular, the limit holds in the case where the forces are Lipschitz: $ \nabla W, \; \nabla V \in W^{1,\infty}$.

However, the Vlasov-Poisson system is not included in this setting, due to the singularity of the Coulomb kernel.
Identifying the torus $\TT^d$ with $\left [- \frac12, \frac12 \right]^d$, with appropriate identifications of the boundary, we note the following properties of the Coulomb kernel $K$. $K \in C^\infty(\TT^d \setminus \{0\})$ is smooth function apart from a point singularity at the origin.
In a neighbourhood of the origin, $K$ can be written in the form
\be \label{def:Coulomb}
K(x) = C_d\frac{x}{|x|^d} + K_0(x), \qquad K_0 \in C^\infty .
\ee
The kernel therefore has a strong singularity of the form $K \sim |x|^{-(d-1)}$.

Forces with a point singularity are of interest in physical applications, since this class includes inverse power laws.
From here on, we discuss forces satisfying bounds of the following form: for some $\beta \in (0, d-1]$,
\be \label{power-law-force}
 \frac{|\nabla W(x)|}{|x|^\beta}\leq C, \frac{|\nabla^2 W(x)|}{|x|^{\beta+1}}\leq C \quad \, \text{ for all } x \in \RR^d \setminus \{0\} .
\ee
Note that the Vlasov-Poisson case corresponds to $\beta = d-1$.

Several works have studied the mean field limit problem for singular forces of the form \eqref{power-law-force} by considering a regularisation of the limit. The singular force $\nabla W$ is replaced by a smooth approximation $\nabla W_r$ such that $\lim_{r \to 0} \nabla W_r = \nabla W$. 
Then, the limits as $N$ tends to infinity and as $r$ tends to zero are taken simultaneously. In this way, one derives the Vlasov equation with singular force in the limit from a sequence of regularised particle systems.
In this formulation, the goal is to optimise the regime $r = r(N)$ for which this limit is valid. That is, $r$ should be as small as possible, so that the regularised particle systems are close to the original particle system with singular interaction.

Hauray and Jabin \cite{Hauray-Jabin} considered the case $\beta < d-1$.
The force is regularised by truncation at a certain distance from the singularity. In this case the regularisation parameter $r(N)$ represents the order of this truncation distance. If $r(N)$ tends to zero sufficiently slowly as $N$ tends to infinity, they prove that the regularised mean field limit holds for a large set of initial configurations. For `weakly singular' forces with $\beta < 1$, in \cite{Hauray-Jabin07, Hauray-Jabin} they also prove the mean field limit without truncation.

For Coulomb interactions, the results available depend on the dimension of the problem. 
In one dimension, the interaction force is less singular. As a consequence, the mean field limit holds, as proved by Hauray \cite{Hauray14}.
The corresponding result for the VPME system was proved by Han-Kwan and the second author in \cite{IHK1}.

In higher dimensions, the Coulomb force is of the form \eqref{power-law-force}. It has a strong singularity corresponding to the endpoint case $\beta = d-1$ not covered by the results of Hauray and Jabin \cite{Hauray-Jabin}. Regularised approaches were considered by Lazarovici \cite{Lazarovici} and Lazarovici and Pickl \cite{Lazarovici-Pickl}. By a truncation method, Lazarovici and Pickl prove a regularised mean field limit for the Vlasov-Poisson system, for a truncation radius of order $r(N) \sim N^{- 1/d + \eta}$ for any $\eta > 0$. To put this in context, note that $N^{- 1/d}$ is the order of separation of particles in $x$ if their spatial distribution is close to uniform.

In a recent breakthrough \cite{Serfaty}, Serfaty introduced a modulated energy method to prove the validity of the mean-field limit for systems of points evolving along the gradient flow of their interaction energy when the interaction is the Coulomb potential or a super-coulombic Riesz potential, in arbitrary dimension. In the appendix (in collaboration with Duerinckx), they adapt this method to prove the mean-field convergence of the solutions to Newton's law with Coulomb interaction in the monokinetic case to solutions of an Euler-Poisson type system.

For the VPME system, a regularised mean field limit was considered by the authors in \cite{GPI20}.
The regularisation used is a regularisation by convolution, similar to the setting of Lazarovici \cite{Lazarovici} that we describe below in Subsection~\ref{sec:VPME-MFL}.
With this regularisation, the resulting microscopic system represents a system of interacting extended charges, where the parameter $r$ gives the order of the radius of the charges. Lazarovici \cite{Lazarovici} derived the Vlasov-Poisson system from a system of extended electrons for $r(N) \geq C N^{-\frac{1}{d(d+2)} + \eta}$ for some $\eta > 0$.
In \cite{GPI20}, the authors proved a similar derivation for the VPME system from a system of extended ions, for the same range of $r$.
We present this result below in Subsection~\ref{sec:VPME-MFL}. 
To our knowledge, this is the first derivation of the VPME system from a particle system in three dimensions.

\subsubsection{Mean Field Limits for VPME} \label{sec:VPME-MFL}

\begin{figure}[h]
\centering
\includegraphics[width=0.4\textwidth]{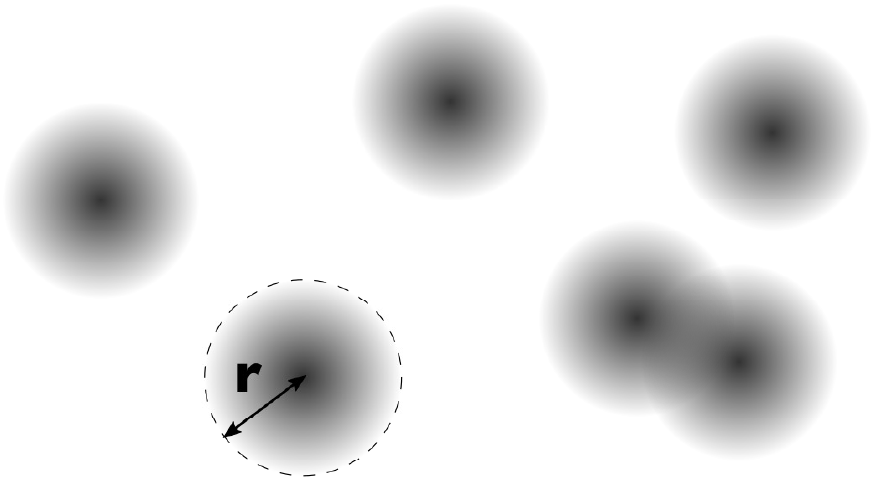}
\caption{A system of extended charges. Here $\chi$ is supported in the unit ball and thus $r$ represents the radius of each charge.}
\end{figure}

For the VPME system, the mean field limit was proved in the one-dimensional setting in \cite{IHK1}.
In the article \cite{GPI20}, we considered the problem in higher dimensions $d=2,3$, deriving the VPME system from a particle system.
The microscopic system is regularised with the regularisation used by Lazarovici \cite{Lazarovici} for the Vlasov-Poisson system.
It consists of a system of `extended ions': instead of representing the ions as point charges, we consider charges of shape $\chi$ for some non-negative, radially symmetric function $\chi \in C^\infty_c(\RR^d)$ with unit mass. The charges are rescaled as follows: for $r > 0$, let
\be \label{def:chi-r}
\chi_r(x) : = r^{-d} \chi \left ( \frac{x}{r} \right ) .
\ee

The extended ions interact with a background of thermalised electrons, leading to the following system of ODEs:
\be \label{eq:ODE-VPME-reg}
\begin{cases}
\dot{X}_i = V_i \\ \ds
\dot{V}_i = - \chi_r \ast \nabla_x U_r(X_i)  , \\ \ds
\Delta U_r = e^{U_r} - \frac{1}{N} \sum_{i=1}^N \chi_r (X_i) .
\end{cases}
\ee

We are able to derive the VPME system \eqref{eq:VPME} from this regularised system, under a condition on the initial data that is satisfied with high probability for $r(N) \geq C N^{-\frac{1}{d(d+2)} + \eta}$. This matches the rate found in Lazarovici's result for the Vlasov-Poisson system. 

\begin{thm}[Regularised mean field limit] \label{thm:MFL-summary}
Let $d=2,3$, and let $f_0 \in L^1 \cap L^\infty(\TT^d \times \RR^d)$ be compactly supported. Let $f$ denote the unique bounded density solution of the VPME system \eqref{eq:VPME} with initial datum $f_0$. Fix $T_* > 0$.

Assume that $r = r(N)$ and the initial configurations for \eqref{eq:ODE-VPME-reg} are chosen such that the corresponding empirical measures satisfy, for some sufficiently large constant $C > 0$, depending on $T_*$ and the support of $f_0$,
\be \label{config-rate}
\limsup_{N \to \infty} \frac{W_2^2(f_0, \mu^N_r(0)) }{r^{d + 2 + C |\log{r}|^{-1/2}}} < 1.
\ee
Then the empirical measure $\mu^N_r$ associated to the particle system dynamics starting from this configuration converges to $f$:
\be \label{thm-statement-MFL-conv}
\lim_{N \to \infty} \sup_{t \in [0,T_*]} W_2(f(t), \mu^N_r(t)) = 0 .
\ee 
In particular, choose $r(N) = N^{-\gamma}$ for some $\gamma < \frac{1}{d(d+2)}$.
For each $N$, let the initial configurations for the regularised $N$-particle system \eqref{eq:ODE-VPME-reg} be chosen by taking $N$ independent samples from $f_0$. Then \eqref{thm-statement-MFL-conv} holds with probability one.
\end{thm}

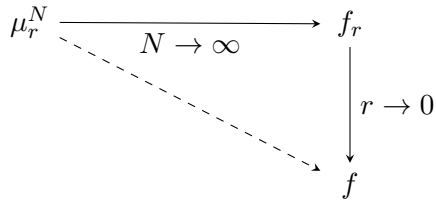
\begin{figure}[ht]
\centering
\begin{tikzpicture}
  \matrix (m) [matrix of math nodes,row sep=4em,column sep=9em,minimum width=2em]
  {
     \mu^N_{r} & f_r \\
       & f \\};
  \path[-stealth]
    (m-1-1)             edge node [below] {$N \rightarrow \infty$} (m-1-2)
                edge [dashed] (m-2-2)
    (m-1-2) edge node [right] {$r \rightarrow 0$} (m-2-2);
\end{tikzpicture}
\caption{Strategy for the proof of Theorem~\ref{thm:MFL-summary}.}
 \label{fig:strat-MFL}
\end{figure}

This theorem is proved by introducing a regularised version of the VPME system:
\be \label{eq:VPME-reg}
\begin{cases}
\partial_t f_r + v \cdot \nabla_x f_r + E_r \cdot \nabla_v f_r = 0, \\
E = - \chi_r \ast_x \nabla_x U, \; \Delta U = e^U - \chi_r \ast_x \rho_f , \\ \ds
f_r \vert_{t=0} = f_0, \;  \int_{\TT^d \times \RR^d} f_0 (x,v) \di x \di v = 1.
\end{cases}
\ee
The solution $f_r$ of this system is used as an intermediate step between the particle system and the VPME system, as illustrated in Figure~\ref{fig:strat-MFL}.

The proof proceeds as follows:
\begin{itemize}
\item We estimate the discrepancy between $\mu^N_r$ and $f_r$, and that between $f_r$ and $f$, in a Wasserstein distance. This uses similar techniques to the stability estimate discussed in Subsection~\ref{sec:strategy-QN}.
\item This estimate is carefully quantified and the regularisation parameter $r$ is allowed to depend on $N$. This allows us to identify a relationship between $r$ and $N$ such that $\mu^N_r$ converges to $f$ for almost all initial data drawn as $N$ independent samples from $f_0$.
\end{itemize}

\subsection{Derivation of Kinetic Euler Systems}

The kinetic Euler systems \eqref{eq:KE-iso}, \eqref{eq:KE-inc} can be derived from particle systems, by using a modified scaling instead of the mean field scaling.
In the articles \cite{IGP1, GPI20} we consider an approach based on a combined mean field and quasineutral limit.
In terms of the scaling $\alpha(N)$, this means that we write $\alpha = (N \e^2)^{-1}$, and then consider allowing $\e$ to depend on $N$. 
We then seek a rate of decay of $\e(N)$ to zero as $N$ tends to infinity for which it possible to take the mean field and quasineutral limits simultaneously.

Due to the challenges involved in the mean field limit for Vlasov-Poisson system, as discussed above, we again use the extended charges model.
For the KIsE system we therefore work with the following microscopic system:
\be \label{eq:ODE-KISE-reg-summary}
 \left \{
\begin{array}{l}
\dot{X}_i = V_i \\ \ds
\dot{V}_i = - \chi_r \ast \nabla_x U(X_i)  , \\ \ds
\e^2 \Delta U = e^{U} - \frac{1}{N} \sum_{i=1}^N \chi_r(x - X_i) .
\end{array}
\right.
\ee
In \cite{GPI20}, we prove the following result.

\begin{thm}[From extended ions to kinetic isothermal Euler] \label{thm:MFQN-KIsE-summary}
Let $d=2$ or $3$, and let $f_\e(0), g_\e(0)$ and $g(0)$ satisfy the assumptions of Theorem~\ref{thm:quasi-summary}. 
Let $T_* > 0$ be the maximal time of convergence from Theorem~\ref{thm:quasi-summary} and let $g$ denote the solution of the KIsE system \eqref{eq:KE-iso} with initial data $g(0)$ on the time interval $[0,T_*]$ appearing in the conclusion of Theorem~\ref{thm:quasi-summary}. 

Let $r = r(N)$ be of the form
\be
r(N) = c N^{-\frac{1}{d(d+2)} + \eta}, \quad \text{for some} \; \; \eta > 0, \; c > 0.
\ee
There exists a constant $C$, depending on $d$, $\eta$, $c$ and $\{ f_\e(0) \}_\e$, such that the following holds.

Let $\e = \e(N)$ satisfy
\be
\e(N) \geq \frac{C}{\sqrt{\log \log \log N}} , \qquad \lim_{N \to \infty} \e(N) = 0.
\ee
For each $N$, let the initial conditions for the regularised and scaled $N$-particle ODE system \eqref{eq:ODE-KISE-reg-summary} be chosen randomly with law $f_{\e(N)}(0)^{\otimes N}$. Let $\mu^N_{\e,r}(t)$ denote the empirical measure associated to the solution of \eqref{eq:ODE-KISE-reg-summary}.

Then, with probability one,
\be
\lim_{N \to \infty} \sup_{t \in [0,T_*]} W_1\left (\mu^N_{\e,r}(t), g(t) \right) = 0.
\ee

\end{thm}

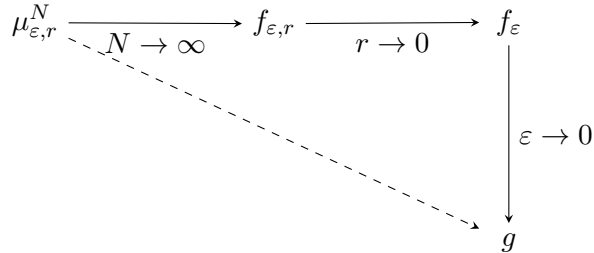
\begin{figure}[ht] 
\centering
\begin{tikzpicture}
  \matrix (m) [matrix of math nodes,row sep=6em,column sep=6em,minimum width=2em]
  {
     \mu^N_{\e,r} & f_{\e,r} & f_{\e} \\
      & & g \\};
  \path[-stealth]
    (m-1-1)             edge node [below] {$N \rightarrow \infty$} (m-1-2)
                edge [dashed] (m-2-3)
    (m-1-2)             edge node [below] {$r \rightarrow 0$} (m-1-3)
    (m-1-3) edge node [right] {$\e \rightarrow 0$} (m-2-3);
\end{tikzpicture}
\caption{Strategy for the proof of Theorem~\ref{thm:MFQN-KIsE-summary}.}
 \label{fig:strat-MFQN}
\end{figure}

This theorem is proved using the strategy illustrated in Figure~\ref{fig:strat-MFQN}. Here $f_{\e,r}$ denotes the solution of a version of the regularised VPME system \eqref{eq:VPME-reg} with quasineutral scaling.

The proof proceeds as follows:
\begin{itemize}
\item As in the proof of Theorem~\ref{thm:MFL-summary}, we estimate the Wasserstein distance between $\mu^N_{\e,r}$ and $f_{\e,r}$ and between $f_{\e,r}$ and $f_\e$.
\item We carefully quantify these estimates in terms of all three parameters $N$, $r$ and here also $\e$.
\item For the convergence of $f_\e$ to $g$, we appeal to Theorem~\ref{thm:quasi-summary}.
\item Using this, we are able to identify a dependence $r = r(N)$ and $\e = \e(N)$ of the parameters on the number of particles, and a relation between $r$ and $\e$, so that the convergence from the particle system to the KIsE system holds for almost all initial data drawn as independent samples from $f_{e}(0)$.
\end{itemize}

\bibliography{proc-MFQN-bib}

\begin{thebibliography}{10}

\bibitem{Baradat}
A.~Baradat.
\newblock Nonlinear instability in {Vlasov} type equations around rough
  velocity profiles.
\newblock {\em Annales de l'Institut Henri Poincar{\'{e}} C, Analyse non
  lin{\'{e}}aire}, 37(3):489 -- 547, 2020.

\bibitem{Bardos}
C.~Bardos.
\newblock {About a Variant of the 1d Vlasov equation, dubbed
  ``Vlasov-Dirac-Benney equation''}.
\newblock In {\em S{\'{e}}minaire Laurent Schwartz - {\'{E}}quations aux
  d{\'{e}}riv{\'{e}}es partielles et applications. Ann{\'{e}}e 2012-2013.},
  S{\'{e}}min. {\'{E}}qu. D{\'{e}}riv. Partielles, pages 1--21. {\'{E}}cole
  Polytechnique, Centre de Math{\'{e}}matiques, Palaiseau, 2014.

\bibitem{Bardos-Besse}
C.~Bardos and N.~Besse.
\newblock {The Cauchy problem for the Vlasov-Dirac-Benney equation and related
  issues in fluid mechanics and semi-classical limits}.
\newblock {\em Kinet. Relat. Models}, 6(4):893--917, 2013.

\bibitem{Bardos-Besse2015}
C.~Bardos and N.~Besse.
\newblock Hamiltonian structure, fluid representation and stability for the
  {V}lasov-{D}irac-{B}enney equation.
\newblock In {\em Hamiltonian partial differential equations and applications},
  volume~75 of {\em Fields Inst. Commun.}, pages 1--30. Fields Inst. Res. Math.
  Sci., Toronto, ON, 2015.

\bibitem{Bardos-BesseSC}
C.~Bardos and N.~Besse.
\newblock Semi-classical limit of an infinite dimensional system of nonlinear
  {S}chr\"{o}dinger equations.
\newblock {\em Bull. Inst. Math. Acad. Sin. (N.S.)}, 11(1):43--61, 2016.

\bibitem{BGNS18}
C.~Bardos, F.~Golse, T.~T. Nguyen, and R.~Sentis.
\newblock The {M}axwell-{B}oltzmann approximation for ion kinetic modeling.
\newblock {\em Phys. D}, 376/377:94--107, 2018.

\bibitem{Bardos-Nouri}
C.~Bardos and A.~Nouri.
\newblock {A Vlasov equation with Dirac potential used in fusion plasmas}.
\newblock {\em J. Math. Phys.}, 53(11):115621, 2012.

\bibitem{Bellan}
P.~M. Bellan.
\newblock {\em Fundamentals of Plasma Physics}.
\newblock Cambridge University Press, 2008.

\bibitem{BNR}
H.~L. Berk, C.~E. Nielsen, and K.~V. Roberts.
\newblock Phase space hydrodynamics of equivalent nonlinear systems:
  Experimental and computational observations.
\newblock {\em The Physics of Fluids}, 13(4):980--995, 1970.

\bibitem{BPLT1991}
G.~Bonhomme, T.~Pierre, G.~Leclert, and J.~Trulsen.
\newblock Ion phase space vortices in ion beam-plasma systems and their
  relation with the ion acoustic instability: numerical and experimental
  results.
\newblock {\em Plasma Physics and Controlled Fusion}, 33(5):507--520, may 1991.

\bibitem{BFJJ}
M.~Bossy, J.~Fontbona, P.-E. Jabin, and J.-F. Jabir.
\newblock {Local existence of analytical solutions to an incompressible
  Lagrangian stochastic model in a periodic domain}.
\newblock {\em Comm. Partial Differential Equations}, 38(7):1141--1182, 2013.

\bibitem{Bouchut}
F.~Bouchut.
\newblock {Global weak solution of the Vlasov-Poisson system for small
  electrons mass}.
\newblock {\em Comm. Partial Differential Equations}, 16(8-9):1337--1365, 1991.

\bibitem{Bouchut-Dolbeault95}
F.~Bouchut and J.~Dolbeault.
\newblock On long time asymptotics of the {Vlasov-Fokker-Planck} equation and
  of the {Vlasov-Poisson-Fokker-Planck} system with {Coulombic} and {Newtonian}
  potentials.
\newblock {\em Differential Integral Equations}, 8(3):487--514, 1995.

\bibitem{Braun-Hepp}
W.~Braun and K.~Hepp.
\newblock {The Vlasov dynamics and its fluctuations in the {$1/N$} limit of
  interacting classical particles}.
\newblock {\em Comm. Math. Phys.}, 56(2):101--113, 1977.

\bibitem{Brenier1989}
Y.~Brenier.
\newblock Une formulation de type {Vlassov--Poisson} pour les {\'{e}}quations
  d'{E}uler des fluides parfaits incompressibles.
\newblock [Rapport de recherche] RR-1070, INRIA, 1989.

\bibitem{Brenier1999}
Y.~Brenier.
\newblock Minimal geodesics on groups of volume-preserving maps and generalized
  solutions of the euler equations.
\newblock {\em Comm. Pure Appl. Math.}, 52(4):411--452, 1999.

\bibitem{Brenier2000}
Y.~Brenier.
\newblock Convergence of the {Vlasov--Poisson} system to the incompressible
  {Euler} equations.
\newblock {\em Comm. Partial Differential Equations}, 25(3-4):737--754, 2000.

\bibitem{Brenier-Grenier94}
Y.~Brenier and E.~Grenier.
\newblock {Limite singuli{\`{e}}re du syst{\`{e}}me de Vlasov-Poisson dans le
  r{\'{e}}gime de quasi neutralit{\'{e}} : le cas ind{\'{e}}pendant du temps}.
\newblock {\em C. R. Acad. Sci. Paris S{\'{e}}r. I Math.}, 318(2):121--124,
  1994.

\bibitem{Carles-Nouri}
R.~Carles and A.~Nouri.
\newblock Monokinetic solutions to a singular {V}lasov equation from a
  semiclassical perspective.
\newblock {\em Asymptot. Anal.}, 102(1-2):99--117, 2017.

\bibitem{Chen}
F.~F. Chen.
\newblock {\em Introduction to Plasma Physics and Controlled Fusion}.
\newblock Springer International Publishing, 3 edition, 2016.

\bibitem{Dobrushin}
R.~L. Dobrushin.
\newblock {Vlasov Equations}.
\newblock {\em Funktsional. Anal. i Prilozhen.}, 13(2):48--58, 1979.

\bibitem{Ferriere}
G.~Ferriere.
\newblock Convergence rate in {Wasserstein} distance and semiclassical limit
  for the defocusing logarithmic schr{\"{o}}dinger equation.
\newblock Preprint, arXiv:1903.04309, 2019.

\bibitem{Golse}
F.~Golse.
\newblock {On the dynamics of large particle systems in the mean field limit}.
\newblock In A.~Muntean, J.~Rademacher, and A.~Zagaris, editors, {\em
  Macroscopic and Large Scale Phenomena: Coarse Graining, Mean Field Limits and
  Ergodicity}, volume~3 of {\em Lect. Notes Appl. Math. Mech.}, pages 1--144.
  Springer, 2016.

\bibitem{Golse-SR2003}
F.~Golse and L.~Saint-Raymond.
\newblock The {V}lasov-{P}oisson system with strong magnetic field in
  quasineutral regime.
\newblock {\em Math. Models Methods Appl. Sci.}, 13(5):661--714, 2003.

\bibitem{Grenier95}
E.~Grenier.
\newblock {Defect measures of the Vlasov-Poisson system in the quasineutral
  regime}.
\newblock {\em Comm. Partial Differential Equations}, 20(7-8):1189--1215, 1995.

\bibitem{Grenier96}
E.~Grenier.
\newblock {Oscillations in quasineutral plasmas}.
\newblock {\em Comm. Partial Differential Equations}, 21(3-4):363--394, 1996.

\bibitem{Grenier99}
E.~Grenier.
\newblock Limite quasineutre en dimension 1.
\newblock In {\em Journ\'{e}es ``\'{E}quations aux {D}\'{e}riv\'{e}es
  {P}artielles'' ({S}aint-{J}ean-de-{M}onts, 1999)}, pages Exp. No. II, 8.
  Univ. Nantes, Nantes, 1999.

\bibitem{IGP-WP}
M.~Griffin-Pickering and M.~Iacobelli.
\newblock Global well-posedness for the {Vlasov-Poisson} system with massless
  electrons in the 3-dimensional torus.
\newblock arXiv:1810.06928.

\bibitem{IGP1}
M.~Griffin-Pickering and M.~Iacobelli.
\newblock A mean field approach to the quasi-neutral limit for the
  {Vlasov--Poisson} equation.
\newblock {\em SIAM J. Math. Anal.}, 50(5):5502--5536, 2018.

\bibitem{GPI20}
M.~Griffin-Pickering and M.~Iacobelli.
\newblock Singular limits for plasmas with thermalised electrons.
\newblock {\em J. Math. Pures Appl.}, 135:199 -- 255, 2020.

\bibitem{GPP}
A.~Gurevich, L.~Pariiskaya, and L.~Pitaevskii.
\newblock Self-similar motion of rarefied plasma.
\newblock {\em Soviet Phys. JETP}, 22(2):449--454, 1966.

\bibitem{GPP2}
A.~Gurevich, L.~Pariiskaya, and L.~Pitaevskii.
\newblock Self-similar motion of a low-density plasma. {II}.
\newblock {\em Soviet Phys. JETP}, 27(3):476--482, 1968.

\bibitem{Gurevich-Pitaevsky75}
A.~V. Gurevich and L.~P. Pitaevsky.
\newblock Non-linear dynamics of a rarefied ionized gas.
\newblock {\em Progress in Aerospace Sciences}, 16(3):227 -- 272, 1975.

\bibitem{Han-Kwan2011}
D.~Han-Kwan.
\newblock Quasineutral limit of the {Vlasov--Poisson} system with massless
  electrons.
\newblock {\em Comm. Partial Differential Equations}, 36(8):1385--1425, 2011.

\bibitem{Han-Kwan-Hauray}
D.~Han-Kwan and M.~Hauray.
\newblock Stability issues in the quasineutral limit of the one-dimensional
  {Vlasov-Poisson} equation.
\newblock {\em Comm. Math. Phys.}, 334(2):1101--1152, 2015.

\bibitem{IHK2}
D.~Han-Kwan and M.~Iacobelli.
\newblock Quasineutral limit for {Vlasov-Poisson} via {Wasserstein} stability
  estimates in higher dimension.
\newblock {\em J. Differential Equations}, 263(1):1--25, 7 2017.

\bibitem{IHK1}
D.~Han-Kwan and M.~Iacobelli.
\newblock {The quasineutral limit of the Vlasov-Poisson equation in Wasserstein
  metric}.
\newblock {\em Commun. Math. Sci.}, 15(2):481--509, 2 2017.

\bibitem{IHK-SCMFL}
D.~Han-Kwan and M.~Iacobelli.
\newblock {From Newton's second law to Euler's equations of perfect fluids}.
\newblock Preprint, arXiv:2006.14924, 2020.

\bibitem{Han-Kwan-Nguyen}
D.~Han-Kwan and T.~T. Nguyen.
\newblock Ill-posedness of the hydrostatic {E}uler and singular {V}lasov
  equations.
\newblock {\em Arch. Ration. Mech. Anal.}, 221(3):1317--1344, 2016.

\bibitem{Han-Kwan-Rousset}
D.~Han-Kwan and F.~Rousset.
\newblock {Quasineutral limit for Vlasov-Poisson with Penrose stable data}.
\newblock {\em Ann. Sci. {\'{E}}c. Norm. Sup{\'{e}}r. (4)}, 49(6):1445--1495,
  2016.

\bibitem{Hauray14}
M.~Hauray.
\newblock Mean field limit for the one dimensional {V}lasov-{P}oisson equation.
\newblock In {\em S\'{e}minaire {L}aurent {S}chwartz---\'{E}quations aux
  d\'{e}riv\'{e}es partielles et applications. {A}nn\'{e}e 2012--2013}, Exp.
  No. XXI, S\'{e}min. \'{E}qu. D\'{e}riv. Partielles. \'{E}cole Polytech.,
  Palaiseau, 2014.

\bibitem{Hauray-Jabin07}
M.~Hauray and P.-E. Jabin.
\newblock {$N$}-particles approximation of the {V}lasov equations with singular
  potential.
\newblock {\em Arch. Ration. Mech. Anal.}, 183(3):489--524, 2007.

\bibitem{Hauray-Jabin}
M.~Hauray and P.-E. Jabin.
\newblock Particle approximation of {Vlasov} equations with singular forces:
  propagation of chaos.
\newblock {\em Ann. Sci. {\'{E}}c. Norm. Sup{\'{e}}r. (4)}, 48(4):891--940,
  2015.

\bibitem{Herda16}
M.~Herda.
\newblock On massless electron limit for a multispecies kinetic system with
  external magnetic field.
\newblock {\em Journal of Differential Equations}, 260(11):7861 -- 7891, 2016.

\bibitem{Jabin-Nouri}
P.~Jabin and A.~Nouri.
\newblock Analytic solutions to a strongly nonlinear {Vlasov} equation.
\newblock {\em C.R. Acad. Sci. Paris, S{\'{e}}r. 1}, 349:541--546, 2011.

\bibitem{Jabin_MFLreview}
P.-E. Jabin.
\newblock A review of the mean field limits for vlasov equations.
\newblock {\em Kinet. Relat. Models}, 7(4):661, 2014.

\bibitem{Lazarovici}
D.~Lazarovici.
\newblock The {Vlasov-Poisson} dynamics as the mean field limit of extended
  charges.
\newblock {\em Comm. Math. Phys.}, 347(1):271--289, 2016.

\bibitem{Lazarovici-Pickl}
D.~Lazarovici and P.~Pickl.
\newblock {A mean field limit for the Vlasov-Poisson system}.
\newblock {\em Arch. Ration. Mech. Anal.}, 225(3):1201--1231, 2017.

\bibitem{Lions-Perthame}
P.~L. Lions and B.~Perthame.
\newblock {Propagation of moments and regularity for the 3-dimensional
  Vlasov-Poisson system}.
\newblock {\em Invent. Math.}, 105(2):415--430, 1991.

\bibitem{Loeper}
G.~Loeper.
\newblock {Uniqueness of the solution to the Vlasov-Poisson system with bounded
  density}.
\newblock {\em J. Math. Pures Appl. (9)}, 86(1):68--79, 2006.

\bibitem{Masmoudi2001}
N.~Masmoudi.
\newblock From {Vlasov--Poisson} system to the incompressible {Euler} system.
\newblock {\em Comm. Partial Differential Equations}, 26(9-10), 2001.

\bibitem{Mason71}
R.~J. Mason.
\newblock Computer simulation of ion-acoustic shocks. {T}he diaphragm problem.
\newblock {\em The Physics of Fluids}, 14(9):1943--1958, 1971.

\bibitem{Medvedev2011}
Y.~V. Medvedev.
\newblock Ion front in an expanding collisionless plasma.
\newblock {\em Plasma Physics and Controlled Fusion}, 53(12):125007, nov 2011.

\bibitem{Mouhot-Villani}
C.~Mouhot and C.~Villani.
\newblock On {L}andau damping.
\newblock {\em Acta Math.}, 207(1):29--201, 2011.

\bibitem{Neunzert-Wick}
H.~Neunzert and J.~Wick.
\newblock {Die Approximation der L{\"{o}}sung von
  Integro-Differentialgleichungen durch endliche Punktmengen}.
\newblock In {\em Numerische Behandlung nichtlinearer Integrodifferential-und
  Differentialgleichungen.}, volume 395 of {\em Lecture Notes in Math.}, pages
  275--290, Berlin, Heidelberg, 1974. Springer.

\bibitem{Penrose}
O.~Penrose.
\newblock {Electrostatic Instabilities of a Uniform Non-Maxwellian Plasma}.
\newblock {\em Phys. Fluids}, 3(2):258--265, 1960.

\bibitem{Pfaffelmoser}
K.~Pfaffelmoser.
\newblock {Global classical solutions of the Vlasov-Poisson system in three
  dimensions for general initial data}.
\newblock {\em J. Differential Equations}, 95(2):281--303, 1992.

\bibitem{SCM}
P.~Sakanaka, C.~Chu, and T.~Marshall.
\newblock Formation of ion-acoustic collisionless shocks.
\newblock {\em The Physics of Fluids}, 14(611), 1971.

\bibitem{Serfaty}
S.~Serfaty.
\newblock Mean field limit for {C}oulomb-type flows.
\newblock {\em Duke Math. J.}, 169(15):2887--2935, 2020.
\newblock Appendix with M. Duerinckx.

\bibitem{Ukai-Okabe}
S.~Ukai and T.~Okabe.
\newblock On classical solutions in the large in time of two-dimensional
  {Vlasov's} equation.
\newblock {\em Osaka J. Math.}, 15(2):245--261, 1978.

\bibitem{Zakharov}
V.~E. Zakharov.
\newblock Benney equations and quasiclassical approximation in the inverse
  problem method.
\newblock {\em Funktsional. Anal. i Prilozhen.}, 14(2):15--24, 1980.

\end{thebibliography}
\bibliographystyle{abbrv}

\end{document}